\documentclass[95pt]{article}
	
\usepackage{verbatim}
\usepackage{hyperref}
\usepackage{amsmath}
\usepackage{bbm}
\usepackage{amsfonts}
\usepackage{amsthm}
\newtheorem*{thm}{Theorem}
\newtheorem{rem1}{Remark}[section]
\newtheorem{lem1}{Lemma}[section]
\newtheorem{cor1}{Corollary}[section]

\newtheorem{thm1}{Theorem}[section]

\begin{document}

\title{Scattering-like phenomena of the periodic defocusing NLS equation}

\author{T. Kappeler\footnote{Supported in part by the Swiss National Science Foundation
and the Newton Institute, Cambridge}, 
B. Schaad\footnote{Supported in part by the Swiss National Science Foundation}, 
P. Topalov\footnote{Supported in part by NSF DMS-0901443}}
\maketitle

\begin{abstract} \noindent 
In this paper we prove approximation properties of the solutions of the defoucsing NLS equation on the circle by nearly linear flows. In addition we show that spatially periodic solutions of the defocusing NLS equation evolving in fractional Sobolev spaces $H^s$ with $s\geq 1$ remain bounded for all times. 
\end{abstract}

\section{Introduction}
Consider the initial value problem for the defocusing  non-linear Schr\"odinger (dNLS) equation 
\begin{align}\label{1} 
i \partial_t u=&-\partial_{x}^2 u + 2 |u|^2u \\\label{2} u_{\vert t=0}=&u_0 
\end{align}  
on the circle $\mathbb{T}=\mathbb{R}/\mathbb{Z}$. According to \cite{Bo0}, 
\eqref{1}-\eqref{2} is globally in time well-posed on the Sobolev spaces
$H^s \equiv H^s(\mathbb{T},\mathbb{C})$ with $s\geq 0$ where
\[
H^s:=\big\{ u= \sum_{n\in \mathbb{Z}}\hat u (n) e^{2\pi i nx }\,\big| \|u\|_{H^s}<\infty \big\}, \quad
\|u\|_{H^s}:=\Big(\sum_{n\in \mathbb{Z}}\langle n\rangle^{2s} |\hat u (n)|^2\Big)^{\frac{1}{2}} , 
\]
with ${\hat u}(n)=\int_0^1 u(x) e^{-2\pi i n x}\,dx$ and $\langle n\rangle:=\max(1,|n|)$. 

 For globally well-posed nonlinear evolution equations on the Euclidean space ${\mathbb R}^n$,
an important issue is whether the solutions admit asymptotics as $t\to \infty$ and in particular,
whether they {\em scatter}, i.e., whether there exist solutions of the linearized equation which approximate the solutions
of the non-linear evolution equation as $t\to \infty$. For the initial value problem \eqref{1}-\eqref{2} on the real line,
it is known that  {\em modified} scattering holds true -- see e.g.  \cite{Carles}, \cite{DZ}, \cite{HN}, \cite{Ozawa}. More precisely, in the setup of weighted Sobolev spaces, it can be shown that the solutions asympotically behave like
$e^{ i S(x,t)} w(x,t)$ for $t \to  \infty$ where $w(x,t) = (4\pi i t)^{-1/2} \int_{\mathbb R} e^{i(x-y)^2 /(4t)} v(y) dy$
with appropriately chosen function $v$. The phase factor $e^{i S(x,t)}$, with real valued phase function $S(x,t),$ reflects the long range effects of the cubic non-linearity
and hence is an asymptotic manifestation of the non-linear dynamics --
see \cite{IT} and references therein for further developments in this direction.
In the periodic setting, due to self-interaction one cannot expect that solutions admit asymptotics as $t\to \infty$.
Our aim is to show that nevertheless, they can be approximated in a sense to be made precise by flows of the form
\begin{align}\label{3}
\sum_{n\in \mathbb{Z}}\hat u_0 (n) e^{-i \omega_n t} e^{2\pi i nx }
\end{align}
where the frequencies $\omega_n$, $n\in \mathbb{Z}$, are real valued and independent of $t$,
but might depend on the initial data $u_0$, encoding in this way part of the non-linear dynamics. 
Note that the flow, defined by \eqref{3} is unitary on $H^s$ and thus leaves the  $H^s$-norm invariant. 

\noindent To state our results let us first recall that the  dNLS equation can be written as a Hamiltonian PDE with the 
real subspace 
\[
H^s_r:= \big\{ (u,\bar{u})\,\big|\,u\in H^s\big\}\subseteq H^s \times H^s
\] 
of the complex space $H^s_c:=H^s\times H^a$ as a phase space and Poisson bracket 
\[
\left\{F,G\right\}:=-i \int_0^1 \big(\partial_u F\partial_{\bar{u}}G- \partial_{\bar u} F \partial_u G\big)\,dx
\]
where $F, G: H^s_r\to\mathbb{C}$ are $C^1$-smooth functionals with sufficiently regular $L^2$-gradients.
The $L^2$-gradients $\partial_u F$ and $\partial_{\bar u}F$ are defined in a standard way in terms of the 
$L^2$-gradients of $F$ with respect to the real and the imaginary part of $u$.
The dNLS equation then takes the form 
\[
\partial_t u= -i \partial_{\bar u}\mathcal{H}_{\rm NLS}
\] 
where $\mathcal{H}_{\rm NLS}$ is the dNLS Hamiltonian 
\[
\mathcal{H}_{\rm NLS}:=\int_0^1 \big(\partial_x u\partial_x \bar u+ u^2 \bar u^2 \big)\,dx.
\]
According to \cite{GK}, the dNLS equation is an integrable PDE in the following strong sense: on $H^0_r$ it admits global canonical 
real analytic coordinates, referred to as {\em Birkhoff coordinates},  so that the initial value problem of the  dNLS equation when 
expressed in these coordinates, can be solved by quadrature. 
These coordinates, denoted by $(z_1, z_2)$, $z_1=\big(z_1(n)\big)_{n\in\mathbb{Z}}$, 
$z_2=\big(z_2(n)\big)_{n\in\mathbb{Z}}$,
are defined on a complex neighborhood of $H^0_r$ in $H^0_c$  and take values in a complex neighborhood of the 
real subspace $\mathfrak{h}^0_r$ of $\mathfrak{h}^0\times\mathfrak{h}^0$ where for any $s\in \mathbb{R}$, 
\[
\mathfrak{h}^s_r:= \big\{(z,{\bar z})\,\big|\, z\in \mathfrak{h}^s\big\}\subseteq
\mathfrak{h}^s_c:=\mathfrak{h}^s\times\mathfrak{h}^s,
\]
and 
$\mathfrak{h}^s:=\mathfrak{h}^s(\mathbb{Z}, \mathbb{C})$ is the Hilbert space
\[
\mathfrak{h}^s:=\big\{z=\big(z(n)\big)_{n\in \mathbb{Z}}\subseteq\mathbb{C} \,\big| \, \|z\|_s < \infty\big\}\, ,
\quad
\|z\|_s:= \Big(\sum_{n\in \mathbb{Z}}\langle n\rangle^{2s} |\,z(n)|^2\Big)^{\frac{1}{2}} \, .
\] 
The Poisson structure on $\mathfrak{h}^s_c$ is defined by the condition that $\left\{ z_1(n), z_2(n) \right\}= -i $ whereas all other 
brackets between coordinate functions vanish. The Birkhoff coordinates have the property that when restricted to $H^N_r$, $N\geq 1$,  
they take values in $\mathfrak{h}^N_r$. 
Furthermore, when expressed in these variables on $\mathfrak{h}^1_r$, the dNLS Hamiltonian is a real analytic function of the actions 
$I_n := z_1(n) \bar z_1(n)$, $n\in \mathbb{Z}$, alone. 
The dNLS equation then reads
\[
\partial_t z_1(n)= \left\{z_1(n), \mathcal{H}_{\rm NLS} \right\}= -i \partial_{{\bar z}_1(n)}\mathcal{H}_{\rm NLS}= 
-i \omega_n^{\rm NLS}z_1(n) \quad \forall n\in \mathbb{Z} 
\]
where
\[
\omega_n^{\rm NLS}:= \partial_{I_n}\mathcal{H}_{\rm NLS} \in 
\mathbb{R}\quad\forall n \in \mathbb{Z}, \: \forall (z_1, \bar z_1)\in \mathfrak{h}^1_r
\]
are referred to as the {\em dNLS frequencies}. 
These frequencies can also be viewed as functions of $u$ which by a slight abuse of notation we also denote by $\omega_n^{\rm NLS}$. 
Note that the frequencies are defined on $H^1$ and are constant along solutions  of \eqref{1}-\eqref{2} in $H^1$.
Our first candidate for the approximate solution of the form \eqref{3}  is 
\begin{align}\label{4}
v(x,t) := \sum_{n\in \mathbb{Z}} \hat u_0(n)e^{- i\omega_n^{\rm NLS}t} e^{2\pi inx}.
\end{align}
The following result makes precise, in which sense $v$ approximates $u$.
\begin{thm1}\label{Theorem1.1} 
For any initial data $u_0 \in H^N$ with $N\in \mathbb{Z}_{\geq 1}$,  the difference $u(t)-v(t) $ is a continuous curve
in $H^{N+1}$ and 
\[
\sup_{t\in \mathbb{R}} \|u(t)-v(t)\|_{H^{N+1}}\leq C < \infty
\] 
where $C>0$  can be chosen uniformly for bounded sets of initial data $u_0\in H^{N}$.
\end{thm1}
\begin{rem1}\label{Remark1.1} Note that Theorem \ref{Theorem1.1} in particular implies that  for any $\epsilon > 0$, $M > 0$, and $N \in \mathbb Z_{\ge 1},$ there exists an integer $L= L_{\epsilon, M, N } > 0$ so that 
the solution $u(t)$ of the dNLS equation \eqref{1} corresponding to initial data 
$u_0$ with $\|u_0\|_{H^{N}}\leq M$,  whose Fourier modes vanish up to order $L$, is approximated
within $\epsilon$ by the nearly linear flow $v(t)$,  
$$
\|u(t)-v(t)\|_{H^{N}}\leq \epsilon \quad   \forall \, t \in \mathbb R .
$$
This result confirms observations made in the optical communication literature that the evolution by \eqref{1} of highly oscillating, localized pulses is nearly linear -- see \cite{EZ} and references therein.
\end{rem1}

\noindent As an immediate application of Theorem \ref{Theorem1.1} we obtain the following  result, 
answering a question concerning the boundedness of solutions of \eqref{1} in fractional  Sobolev spaces, raised in 
Dispersive Wiki \cite{DiWi}. Actually, the boundedness can be shown to be uniform in the following sense:
\begin{cor1}\label{Corollary1.2}
For any $s\in \mathbb{R}_{\geq 1}$, the solution $u(x,t)$ of \eqref{1}-\eqref{2} with initial data in  $H^s$ stays bounded for all times, 
\begin{equation}\label{eq:uniform}
\sup_{t\in \mathbb{R}} \|u(t)\|_{H^s}\leq C < \infty
\end{equation}
where $C>0$ can be chosen uniformly on bounded sets of initial data in $H^s.$
\end{cor1}
\noindent We remark that the novelty of Corollary \ref{Corollary1.2} is that \eqref{eq:uniform} holds
in fractional Sobolev spaces. In the case where $s$ is an integer the inequality \eqref{eq:uniform} is well known and can be proved
using the conservation laws of the NLS hierarchy.

\smallskip

\noindent It is also natural  to investigate the linear flow obtained when $\omega_n$, $n\in \mathbb{Z}$, in \eqref{3} are
chosen to be the frequencies $4\pi^2 n^2, \, n\in \mathbb{Z}$, of the linear Schr\"odinger equation on the circle. 
Similarly as in the case of the modified scattering on the line, it turns out to be better to modify the flow of the
linear Schr\"odinger equation by a phase factor. Consider the following modified version of the linear Schr\"odinger equation,
\begin{align}
\label{5} 
i \partial_t w=&  -\partial_x^2 w + P(w)w, \quad P(w):= 4 \int_0^1 |w(x,t)|^2  dx 
 \\
\label{6}
w_{\vert t=0}=& u_0 \, .
\end{align} 
Note that the solution $w(t,x)$ of \eqref{5}-\eqref{6} for $u_0\in H^s$ with $s\geq 0$ leaves the $L^2$-norm invariant, 
\begin{align*} c:= 4 \int_0^1|u_0(x)|^2dx= 4 \int_0^1 |w(x,t)|^2dx\quad \forall t\in \mathbb{R}\end{align*} 
and is given by 
\begin{align}\label{7}
w(x,t)= e^{-ict} \sum_{n\in \mathbb{Z}} \hat u_0(n)e^{-i4\pi^2 n^2t}e^{2\pi inx}.
\end{align}
It is our second candidate for a flow of the form \eqref{3}.
\begin{thm1}\label{Theorem1.3}
For any initial data $u_0\in H^N$ with $N\in \mathbb{Z}_{\geq 2}$, 
\[ 
\|u(t)-w(t)\|_{H^{N+1}} \leq C(1+ |t|) \quad \forall t\in \mathbb{R}.
\]
The constant $C>0$ can be chosen uniformly on bounded sets of initial data $u_0\in H^N$.
\end{thm1}
\begin{rem1}\label{Remark1.2} Note that Theorem \ref{Theorem1.3} in particular implies that  for any $\epsilon > 0$, $M > 0$, $T > 0$, and $N \in \mathbb Z_{\ge 2},$ there exists an integer $L= L_{\epsilon, M, T, N } > 0$ so that 
the solution $u(t)$ of the dNLS equation \eqref{1} corresponding to initial data 
$u_0$ with $\|u_0\|_{H^{N}}\leq M$,  whose Fourier modes vanish up to order $L$, is approximated
on the time interval $[-T, T]$ within $\epsilon$ by the nearly linear flow $w(t)$,
$$
\|u(t)-w(t)\|_{H^{N}}\leq \epsilon \quad   \forall \, t \in [-T, T] .
$$
This result again confirms observations made in the optical communication literature that the evolution by \eqref{1} of highly oscillating, localized pulses is nearly linear -- see Remark \ref{Remark1.1}.
\end{rem1}

\medskip

\noindent {\em Outline of proofs:} The proof of Theorem \ref{Theorem1.1}, given in Section 3, relies on the semi-linear property of the 
Birkhoff coordinates established in \cite{KSTnls}; see Section 3 for a review. 
This property leads to additional features of the difference $u(t)-v(t)$ which are discussed in the same section. 
The proof of Theorem \ref{Theorem1.3}, given also in Section 3, is based on the following asymptotics of the dNLS  frequencies 
$\omega_n^{\rm NLS}$ for $|n|\to \infty$, proved in Section 2 where we also introduce additional notation and review the setup. 
\begin{thm1}\label{Theorem1.4}
For any $u\in H^2$,
\[\omega_n^{\rm NLS}= 4 n^2\pi^2+ 4 \int_0^1|u(x)|^2dx+ O\Big(\frac{1}{n}\Big) \quad \text{as} \: |n|\to \infty\]
uniformly on bounded subsets of $H^2$.
\end{thm1}
\noindent In Theorem \ref{Theorem2.7}, we state and prove a stronger version of Theorem \ref{Theorem1.4},
valid in a complex neighborhood of $H^2$.

\medskip

\noindent {\em Related work:} The question of growth of Sobolev norms of solutions of 
non-linear dispersive evolution equations such as the dNLS equation on the circle has been studied intensively in recent years. 
For solutions of many equations of this type, polynomial bounds in $t$ of Sobolev norms 
of solutions were established -- e.g. \cite{Bo1}, \cite{Bo2}, \cite{St}; in \cite{So} such bounds 
were found for solutions of the defocusing NLS equation on the circle with a non-linear term  of the form $|u|^{2k}u$.
In \cite{CKSTT}, in a different direction of research, solutions of the defocusing NLS equation 
on the $2d$-torus were constructed with the property that corresponding Sobolev 
norms get very large as time evolves. This work led to further refinements and extensions. 
The results presented in this paper continue a direction of research initiated by Erdogan, Tzirakis, 
and Zharnitsky 
(\cite{ETZ1}, \cite{ET1}, \cite{ET2}, \cite{EZ}; cf also \cite{O}) and the authors \cite{KST1}. 
In particular, in \cite{KST1}, results similar to the ones of Theorem \ref{Theorem1.1}, 
Corollary \ref{Corollary1.2} and Theorem \ref{Theorem1.3} have been established for the KdV equation. 
In \cite{EZ}, \cite{ET2}, one finds results on the approximation of solutions of the focusing NLS equation, related to Theorem \ref{Theorem1.3} and in the appendix in \cite{EZ}, a discussion on related results in optical communication with further references -- see also Remark \ref{Remark1.1} above.
Specifically, let us mention in detail the following result of Erdogan and Tzirakis in \cite{ET2},
proved for the focusing NLS equation, which in the notation of this paper reads as follows:

\begin{thm}[\cite{ET2}, Theorem 2]
For any $u_0\in H^s$ with $s>0$, $u-w$ is in $C^0(\mathbb{R}, H^{s+a})$ for any 
$a< \min\big(2s,\frac{1}{2}\big)$ where here $u$ denotes the solution of the focusing NLS equation with initial data $u_0$ and $w$ the approximation corresponding to \eqref{5}.
\end{thm}
\noindent Most likely, the above Theorem also holds in the case of the defocusing NLS equation. Note that in comparison,  Theorem \ref{Theorem1.3} above states that for solutions of the defocusing NLS equation
in $H^N$ with $N \ge 2,$ $a=1$ and $ \|u(t)-w(t)\|_{H^{N+1}}$ is bounded by $C(1+|t|)$. 

\medskip

\noindent {\em Acknowledgment:} It is a pleasure to thank R\'emi Carles and Vedran Sohinger for 
pointing out to us various references concerning the scattering theory of the NLS equation.

\section{Preliminaries}
In this section we introduce additional notation and discuss in more detail the Birkhoff coordinates and the frequencies, 
introduced in Section 1. In addition we prove Theorem \ref{Theorem1.4}. 
We begin by defining the following Fourier transform,
\[
\mathcal{F}:H^0_c\to \mathfrak{h}^0_c, \quad (\varphi_1, \varphi_2) \mapsto 
\left(\left(-\hat\varphi_1(-n)\right)_{n\in\mathbb{Z}},\left(-\hat\varphi_2(n)\right)_{n\in\mathbb{Z}}\right)\,.
\]
Note  that for $(\varphi_1, \varphi_2)\in H^0_r$, one has $\bar\varphi_2= \varphi_1$, implying that 
$\overline{\hat\varphi}_2(n)= \hat\varphi_1(-n)$ for any $n\in \mathbb{Z}$. Hence $\mathcal{F}$ maps $H^0_r$ in 
$\mathfrak{h}^0_r$. Clearly, for any $s\geq 0$, $\mathcal{F}: H^s_r \to \mathfrak{h}^s_r$ is an isometry. In \cite{GK}, 
Birkhoff coordinates $x(n), y(n),\; n\in\mathbb{Z} $ are constructed on $H^0_r$. 
These coordinates are real valued and satisfy $\left\{x(n), y(n)\right\}=-1$ for any $n\in \mathbb{Z}$ whereas all other brackets 
between coordinate functions vanish. For our purposes, in this paper it is convenient to use complex coordinates
\[
z_1(n)=\frac{x(n)-i y(n)}{\sqrt{2}}\quad \text{and} \quad 
z_2(n)= \frac{x(n)+ i y(n)}{\sqrt{2}} \, (= \bar z_1(n)).
\]
Then $\{z_1(n), \bar z_1(n)\}=-i$ and the action variable $I_n$ can be expressed as 
\[
I_n:= \big( x(n)^2+ y(n)^2\big) / 2 = z_1(n)\bar z_1(n)\quad  \forall n\in \mathbb{Z}.
\]
The result on the Birkhoff coordinates in \cite{GK} (see the Overview as well as Theorem 20.2)
then reads as follows:
\begin{thm1}\label{Theorem2.1}
There exists a real analytic diffeomorphism $\Phi: H^0_r\to \mathfrak{h}^0_r$, 
$(u,\bar u)\mapsto (z_1,\bar z_1)$ with the following properties:
\begin{itemize}
\item[(B1)] $\Phi$ is canonical, i.e., preserves the Poisson brackets. 
\item[(B2)] The restriction of $\Phi$ to $H^N_r, \, N\in \mathbb{Z}_{\geq 1}$, gives rise to a map $\Phi: H^N_r\to \mathfrak{h}^N_r$ that is onto and bianalytic. 
\item[(B3)] $\Phi$ defines global Birkhoff coordinates for the dNLS equation on $H^1_r$. More precisely, 
on $\mathfrak{h}^1_r$, the dNLS Hamiltonian $\mathcal{H}_{\rm NLS}\circ \Phi^{-1}$ is a real analytic function of the actions $I_n=z_1(n)\bar z_1(n), \, n\in \mathbb{Z},$ alone.
\item[(B4)] The differential of $\Phi$ at $0$ is the Fourier transform, $d_0\Phi= \mathcal{F}$.
\end{itemize}  
The map $\Phi$ is referred to as Birkhoff map.
\end{thm1}
\begin{rem1}\label{Remark2.2}
The map  $\Phi:H^0_r\to \mathfrak{h}^0_r$ being real analytic means that there exists a neighborhood $W$ of $H^0_r$ in the complex Hilbert space $H^0_c$ $(= H^0\times H^0)$
where the map $\Phi$ is defined and takes values in the complex Hilbert space $\mathfrak{h}^0_c$ $(=\mathfrak{h}^0\times \mathfrak{h}^0), \: \Phi : W\to \mathfrak{h}^0_c,$ so that $\Phi$ is analytic on $W$.  The coordinates defined by $\Phi$ are denoted by $z_1(n),z_2(n),\,n\in \mathbb{Z}$. In general, $z_2(n)\neq \bar z_1(n)$. The neighborhood can be chosen in such a way that $\Phi(W\cap H^N_c)\subseteq\mathfrak{h}^N_c$ for any $N\in \mathbb{Z}_{\geq 1}$.
\end{rem1}
\begin{rem1}\label{Remark2.3}
The claim (B4) of Theorem \ref{Theorem2.1} follows from \cite{GK},  Theorem 17.2. 
\end{rem1}
\noindent Now let us turn to the dNLS frequencies, introduced in Section 1. The starting point of our investigations is a representation of 
the frequencies in terms  of spectral data associated to the Zakharov-Shabat 
operator -- see \cite{GK} for a detailed analysis. For $\varphi\in H^0_r$, let  \[L(\varphi) = i \begin{pmatrix}
1&0\\0&-1
\end{pmatrix}\partial_x + \begin{pmatrix}
0& \varphi_1 \\
\varphi_2 &0
\end{pmatrix}\]
be the Zakharov-Shabat operator. Denote by
$\cdots \leq \lambda^-_{n}\leq \lambda^+_n< \lambda^-_n \leq\lambda_{n+1}^+ \leq \cdots$
the periodic spectrum of $L(\varphi)$ on $[0,2]$ listed with multiplicities in such a way that $\lambda^{\pm}_n=n\pi + o(1)$ as $n\to \infty$. Furthermore, let 
$\tau_n= (\lambda^+_n+ \lambda^-_n) / 2$ and $\gamma_n= \lambda_n^+-\lambda_n^-$. 
Denote by $\Delta(\lambda)$ the discriminant of $L(\varphi)$ and recall that $\Delta^2(\lambda)-4$ 
admits the product representation 
\[
\Delta^2(\lambda)-4= -4 \prod_{n\in \mathbb{Z}}\frac{(\lambda_n^+-\lambda)(\lambda_n^-- \lambda)}{\pi_n^2}
\]
where $\pi_n= n\pi$ for $n\neq 0$ and $\pi_0= 1$. In addition denote by $\psi_n(\lambda),\,n\in \mathbb{Z},$ the entire function of the the form
\begin{align}
\label{2.2} \psi_n(\lambda)= - \frac{2}{\pi_n} \prod_{k\neq n}\frac{\sigma^n_k-\lambda}{\pi_k}
\end{align}
where the normalizing factor $-\frac{2}{\pi_n}$ and the zeros $\sigma^n_k,\,k\neq n$, are determined in such a way that for any $n,m\in \mathbb{Z}$
\begin{align}
\label{2.3} \frac{1}{2\pi}\int_{\Gamma_m} \frac{\psi_n(\lambda)}{\sqrt[c]{\Delta^2(\lambda)-4}}d\lambda=\delta_{nm}.
\end{align}
Here $\Gamma_m$ is the contour around the interval $[\lambda^-_m,\lambda^+_m]$ so that all other periodic eigenvalues of 
$L(\varphi)$ are outside of $\Gamma_m$, $\sqrt[c]{\Delta^2(\lambda)-4}$ denotes the canonical root 
\begin{align}\label{2.4}
\sqrt[c]{\Delta^2(\lambda)-4}=2i \prod_{n\in \mathbb{Z}}\frac{w_k(\lambda)}{\pi_k}, \quad 
w_k(\lambda)=\sqrt[s]{(\lambda_k^+-\lambda)(\lambda_k^-- \lambda)}
\end{align}
and $\sqrt[s]{(\lambda_k^+-\lambda)(\lambda_k^-- \lambda)}$ is the standard root on $\mathbb{C}\setminus  [\lambda^-_k,\lambda^+_k]$ defined by setting
\[
\sqrt[s]{(\lambda_k^+-\lambda)(\lambda_k^-- \lambda)}
= (\tau_k-\lambda)\sqrt[+]{1- \Big(\frac{\gamma_k/2}{\tau_k-\lambda}\Big)^2}
\quad \text{ for }\: |\lambda|\: \text{ large}.
\]
The zeros $\sigma^n_k,\, k\neq n$, of \eqref{2.2} are listed in such a way that $\lambda_k^-\leq \sigma_k^n\leq \lambda_k^+$ for any $k \neq n$. In the construction of the Birkhoff map in \cite{GK}, the complex neighborhood $W$ of $H^0_r$ in $H^0_c$, mentioned in Remark \ref{Remark2.2}, is chosen in such a way that the spectral analysis of $L(\varphi)$ can be extended for $\varphi\in W$ so that $\tau_n, \gamma_n^2,$ and $\sigma_k^n,\,k\neq n,$ are analytic functions on $W$. The following results are established in  \cite{KST2} and \cite{GK}, Lemma 14.12.
\begin{thm1}\label{Theorem2.4}
On $W\cap H^1_c,$ the spectral quantities $\tau_k, \gamma_k^2,$ and  $\sigma_k^n,\, n\neq k$ are real analytic for any $k\in \mathbb{Z}$ and satisfy
(i) $\tau_k= k\pi + \frac{1}{2\pi k} \int_0^1 \varphi_1(x)\varphi_2(x) dx + O \left(\frac{1}{k^2}\right)$;
(ii) $\gamma_k^2= \frac{1}{k^2} \ell^1(k)$; 
(iii) $\sigma^n_k= \tau_k + O\left(\gamma_k^2\right)$.
These estimates hold uniformly in $n$ and on bounded subsets of $W\cap H^1_c$.
\end{thm1}
\noindent In particular, the estimates of Theorem \ref{Theorem2.4} imply 
\begin{align}\label{2.5}
\tau_n^2= n^2\pi^2 + \int_0^1 \varphi_1(x)\varphi_2(x) dx + O\Big( \frac{1}{n}\Big) \\
\label{2.6} \sum_{k}|\gamma_k|^2 = O(1) \quad\text{and} \quad \sum_{k\neq n}\left|(\sigma^n_k- \tau_k)(\sigma^n_k+\tau_k)\right|= O(1)
\end{align} 
uniformly in $n$ and on bounded subsets of $W\cap H^1_c$. 
By \cite{GK} Theorem 5.7, Lemma 5.8, the dNLS frequencies  
$\omega_n^{\rm NLS}= \partial\mathcal{H}_{\rm NLS}/\partial I_n$ are given on $H^1_r$ by the formula 
\begin{align}
\label{2.7}
\omega_n^{\rm NLS}=2 \tau_n^2 + 2n^2\pi^2+ 2\sum_{k\neq n}(\sigma^n_k- \tau_k)(\sigma^n_k+\tau_k)+ \frac{1}{2} \sum_k \gamma_k^2.
 \end{align}
  Note that due to \eqref{2.6}, the infinite sums in \eqref{2.7} are absolutely convergent on $W\cap H^1_c$ uniformly in $n$ and on bounded subsets of $W\cap H^1_c$. 
 As both sides of \eqref{2.7} are analytic on $W\cap H^1_c,$ identity \eqref{2.7} actually holds on W. By  \eqref{2.5} it then follows that  
 \begin{align}\label{2.8}
 \sup_{n\in \mathbb{Z}} |\omega_n(\varphi)|\leq C \quad \mbox{where} \quad  
\omega_n:= \omega_n^{\rm NLS}- 4\pi^2n^2\quad \forall n \in \mathbb{Z}
 \end{align}
and $C$ can be chosen uniformly on bounded subsets of $W\cap H^1_c$. Let $\ell^\infty\equiv \ell^\infty(\mathbb{Z}, \mathbb{C})$. A (possibly non-linear) map between (subsets of) Banach spaces is said to be bounded if the image of any bounded subset in the domain of definition is bounded. The discussion above leads to the following 
 \begin{cor1}\label{Corollary2.5} 
 The map $W\cap H^1_c \to \ell^\infty, \; \varphi \mapsto\big(\omega_n( \varphi)\big)_{n\in \mathbb{Z}}$ is real analytic and bounded. 
 \end{cor1} 
\begin{proof}[Proof of Corollary \ref{Corollary2.5} ] For any $n\in \mathbb{Z}, \: \omega_n: W\cap H^1_c\to \mathbb{C}$ is analytic 
and real valued on $H^1_r$. By \eqref{2.8}, $(\omega_n)_{n\in \mathbb{Z}}\in \ell^\infty$ and the map 
$\varphi\to\big(\omega_n( \varphi)\big)_{n\in \mathbb{Z}}$ is bounded. Hence by \cite{KdV&KAM}, Theorem A.3, the claimed statement
 follows. 
 \end{proof}
\noindent For Section 3, introduce for any $t\in \mathbb{R}$ and $\varphi \in W\cap H^1_c$ 
\begin{align}
\label{2.9} \big(\Omega^t(\varphi)\big)_{n\in \mathbb{Z}} = 
\Big(\big(\Omega^t_1(n,\varphi)\big)_{n\in \mathbb{Z}},
\big(\Omega^t_2(n,\varphi)\big)_{n\in \mathbb{Z}}\Big) 
\end{align}
where \begin{align}\label{2.10}
\Omega^t_1(n,\varphi)= \exp\big(i \omega_n^{\rm NLS}(\varphi) t\big)\quad 
\text{and} \quad \Omega^t_2(n,\varphi)= \exp\big(-i \omega_n^{\rm NLS}(\varphi) t\big).
\end{align}
Note that by Corollary \ref{Corollary2.5}, $\Omega^t(\varphi)\in \ell^\infty_c:= \ell^\infty\times \ell^\infty$ for $\varphi\in W\cap H^1_c$ whereas for $\varphi\in H^1_r,$
$\Omega^t(\varphi)\in \ell^\infty_r:= \{ (z_1,z_2) \in \ell^\infty_c\big| \, \bar z_2(n) = z_1(n)\}$. 
Corollary \ref{Corollary2.5} actually implies the following 
\begin{cor1}
For any $t\in \mathbb{R}$, the map $\Omega^t: W\cap H^1_c$ is real analytic and bounded.
\end{cor1}
\noindent We finish this section with proving Theorem \ref{Theorem1.4}. We show the following slightly stronger version. 
\begin{thm1}\label{Theorem2.7}
Let $W$ be the neighborhood of $H^0_r$ in $H^0_c$ given in Remark \ref{Remark2.2}. 
Then the map 
$W\cap H^2_c \to \ell^\infty, \; \varphi \mapsto\big((1+|n|)\rho_n( \varphi)\big)_{n\in \mathbb{Z}}$ with
\[
\rho_n( \varphi) : = \omega_n^{\rm NLS}(\varphi) - 
\Big( 4 n^2\pi^2+ 4 \int_0^1\varphi_1(x)\varphi_2(x)\,dx\Big)
\]
 is real analytic and bounded.  In particular, 
\[
\omega_n^{\rm NLS}(\varphi) = 4 n^2\pi^2+ 4 \int_0^1\varphi_1(x)\varphi_2(x) dx
+ O\Big(\frac{1}{n}\Big)\quad \text{as}\; |n|\to \infty 
\]
uniformly on bounded subsets of $W\cap H^2_c$.
\end{thm1}
\begin{proof}[Proof]
Our starting point is formula \eqref{2.7}. By Theorem \ref{Theorem2.4} and \eqref{2.5} 
\begin{align*}
\omega_n^{\rm NLS}= & 4 n^2\pi^2+ 2 \int_0^1 \varphi_1(x)\varphi_2(x)dx \\ 
&+ 2 \sum_{|k|\leq |n| / 2}\Big( (\tau_k-\sigma_k^n)(\tau_k+\sigma_k^n)+ \Big(\frac{\gamma_k}{2}\Big)^2\Big)+ 
O\Big(\frac{1}{n}\Big)
\end{align*}
uniformly on bounded subsets of $W\cap H^1_c$.
The claimed asymptotic estimate then follows from the estimate of Lemma \ref{Lemma2.10} below, 
\begin{align}\label{2.11}
\sum_{|k|\leq |n| / 2}\Big( (\tau_k-\sigma_k^n)(\tau_k+\sigma_k^n)+\Big(\frac{\gamma_k}{2}\Big)^2\Big)=
\int_0^1 \varphi_1(x)\varphi_2(x)\,dx +O\Big(\frac{1}{n}\Big)
\end{align}
which holds uniformly on bounded subsets of $W\cap H^2_c.$
As for each $n \in \mathbb Z,$ $\rho_n : W\cap H^2_c \to \mathbb C$ is real analytic, 
one then concludes that the map $W\cap H^2_c \to \ell^\infty$,
$\varphi\mapsto\big((1+|n|)\rho_n( \varphi)\big)_{n\in \mathbb{Z}}$ is real analytic.
\end{proof}
\noindent It remains to show the asymptotic estimate \eqref{2.11}, used in the proof of Theorem \ref{Theorem2.7}. To this end we consider the function $\prod_{|k|\leq |n| / 2}\frac{\sigma_k^n-\lambda}{w_k(\lambda)}$.
Substitute $\lambda=- \frac{1}{z}$ in the latter product and define for $|z|\leq |\tau_n|^{-1}$ with $|n|$ sufficiently large 
\begin{align}
\label{2.11bis} 
g_n(z):= \prod_{|k|\leq |n| / 2} \frac{1+ z\sigma_k^n}{\sqrt[+]{(1+ z \lambda_k^+)(1+ z \lambda_k^-)}}.
\end{align}
In Lemma \ref{Lemma2.9} below we compute its expansion of order $3$ at $z=0$. Combining the identity
$\sum_{k\ne n} (\sigma_k^n-\tau_k)=\tau_n-n\pi$
established in \cite{GK1}, Lemma 5.8, with the asymptotic estimate of 
$\prod_{k\neq n}\frac{\sigma_k^n-\tau_n}{w_k(\lambda)}$ for $|n|\to \infty,$ stated in Lemma \ref{Lemma2.8} below, we will prove \eqref{2.11}. We begin with the estimate of the latter product. Let 
\[
f_n(\mu):=\prod_{k\neq n}\frac{\sigma_k^n-\lambda}{w_k(\lambda)},\quad \mu=\lambda-\tau_n.
\]
\begin{lem1}
\label{Lemma2.8}
On $W\cap H^1_c,\; f_n(0)= 1+ O\left(\frac{1}{n^3}\right)$ uniformly on bounded subsets of $W\cap H^1_c$. 
\end{lem1}
\begin{proof}[Proof of Lemma \ref{Lemma2.8}] By the counting lemma (cf. e.g. \cite{GK}, Lemma 6.4) there exists $n_0\geq 1$ so that for any $|n|>n_0,\, f_n$ is analytic on the disc 
$D_n=\left\{|\lambda-n\pi|< \frac{\pi}{4}\right\}$ and $|\lambda_n^\pm-n\pi|\leq \frac{\pi}{8}$. 
The number $n_0$ can be chosen uniformly on bounded subsets of $W\cap H^1_c.$ By \eqref{2.3} one has 
\[
\frac{1}{2\pi}\int_{\Gamma_n} \frac{\psi_n(\lambda)}{\sqrt[c]{\Delta^2(\lambda)-4}}d\lambda=1.
\]
In case where $\lambda_n^-= \lambda_n^+$ one has by Cauchy's formula
\[1= \frac{1}{2\pi i}\int_{\Gamma_n}\frac{f_n(\lambda-\tau_n)}{\lambda-\tau_n}d\lambda= f_n(0).\]
If $\lambda_n^-\neq \lambda_n^+$ with $|n| >n_0,$ 
argue as in the proof of \cite{GK}, Theorem 13.3 to see that
\begin{align} \label{2.12}
1=\frac{1}{\pi}\int_0^1\Big(f_n\Big(t \frac{\gamma_n}{2}\Big)+f_n\Big(-t\frac{\gamma_n} {2}\Big)\Big)
\frac{dt}{\sqrt[+]{1-t^2}}\,.
\end{align}
Expanding $f_n$ near $\mu=0$, one gets
$f_n(\mu)=f_n(0)+ f'_n(0)\mu+ \mu^2 \int_0^1f_n''(s\mu)(1-s)\,ds$ yielding 
\[
f_n(\mu) + f_n(-\mu)= 2f_n(0)+ \mu^2\int_0^1 \big(f_n''(s\mu)+f_n''(-s\mu)\big)(1-s)\,ds.
\]
When substituted into \eqref{2.12} one obtains 
$1=\frac{2f_n(0)}{\pi} \int_0^1\frac{dt}{\sqrt[+]{1-t^2}} + A$ where
\[
A:= \frac{1}{\pi} \int_0^1 \left(\int_0^1\Big(f_n''\Big(s t \frac{\gamma_n}{2}\Big)+ 
f_n''\Big(-s t \frac{\gamma_n}{2}\Big)\Big)(1-s)\,ds\right)
\frac{t^2\big(\gamma_n /2 \big)^2}{\sqrt[+]{1-t^2}}dt.
\]
As $\int_0^1\frac{dt}{\sqrt[+]{1-t^2}}=\frac{\pi}{2}$ it then follows that
\[
|f_n(0)-1|\leq\Big|\frac{\gamma_n}{2}\Big|^2\max_{-1\leq s\leq 1}\Big|f_n''\Big(s\frac{\gamma_n}{2}\Big)\Big|\,.
\]
By Theorem \ref{Theorem2.4}, $\gamma_n^2= O\big(\frac{1}{n^2}\big)$ uniformly on bounded subsets of $W\cap H^1_c$. As for $|n|> n_0,\, |\lambda_n^\pm-n\pi|\leq \frac{\pi}{8}$ and $f_n(\lambda-\tau_n)$ is analytic in $\lambda$ on the disc $D_n$ one can bound 
$f_n''(\mu) =\partial_\mu^2\left(f_n(\mu)-1\right)$ on $[\lambda_n^-,\lambda_n^+]= \left\{(1-t)\lambda_n^-+ t\lambda_n^+|\,0\leq t\leq 1\right\}$  by Cauchy's estimate 
\[
\max_{-1\leq s\leq 1} \left|f''_n\left(s\frac{\gamma_n}
{2}\right)\right|\leq C \sup_{\lambda \in D_n}\left|f_n(\lambda-\tau_n)-1\right|
\]
where $C>0$ can be chosen uniformly on bounded subsets of $W\cap H^1_c$. It remains to prove that 
\begin{align}\label{2.13}
\sup_{\lambda \in D_n}\left|f_n(\lambda-\tau_n)-1\right|=O\Big(\frac{1}{n}\Big)
\end{align}
uniformly on bounded subsets of $W\cap H^1_c$.
For convenience, let $\alpha_k^n=\sigma_k^n-\tau_k\; (k\neq n)$ and write $f_n(\lambda-\tau_n)$ with $|n|>n_0$ as follows 
\[
f_n(\lambda-\tau_n)= \prod_{k\neq n} \frac{\sigma_k^n-\lambda}{\sqrt[s]{(\tau_k-\lambda)^2-\frac{\gamma_k^2}{4}}}=
\prod_{k\neq n} \frac{1+ \frac{\alpha_k^n}{\tau_k-\lambda}}{\sqrt[+]{1- \big(\frac{\gamma_k / 2}{\tau_k-\lambda}\big)^2}}.
\]
Choose $n_0$ larger if necessary so that for any $|n|>n_0$ and any $\lambda\in \bar D_n$
\[
\Big|\frac{\alpha_k^n}{\tau_k-\lambda}\Big|\leq \frac{1}{2}\quad\forall k\neq n\quad \text{and} \quad
\Big|\frac{\gamma_k}{\tau_k-\lambda}\Big|\leq\frac{1}{2}\quad\forall k\ne n.
\] 
In view of Theorem \ref{Theorem2.4}, $n_0$ can be chosen uniformly on bounded subsets of $W\cap H^1_c$. 
By Theorem \ref{Theorem2.4} and \cite{GK}, Lemma C.2  it follows that 
\begin{align*}
\Big|\prod_{k\neq n}\Big(1+ \frac{\alpha_k^n}{\tau_k-\lambda}\Big)-1\Big|
=\Big|\exp\Big(\sum_{k\neq n}\log\Big(1+\frac{\alpha_k^n}{\tau_k-\lambda}\Big)-1\Big)\Big|
=O\Big(\frac{1}{n}\Big)
\end{align*}
and 
\[
\Big|\prod_{k\neq n}\Big(1-\Big(\frac{\gamma_k / 2}{\tau_k-\lambda}\Big)^2\Big)-1\Big|=  
\Big|\exp\Big(\sum_{k\neq n}\log\Big(1-\Big(\frac{\gamma_k / 2}{\tau_k-\lambda}\Big)^2\Big)-1\Big)\Big|
=O\Big(\frac{1}{n^2}\Big)
\]
uniformly on bounded subsets of $W\cap H^1_c.$ Combining these estimates one concludes that \eqref{2.13} holds. 
\end{proof}
\noindent Let us now turn towards the function $g_n$ introduced in \eqref{2.11bis}. From the asymptotics of 
$\tau_k, \, \gamma_k, $ and $\sigma_k^n$ it follows that for any $|n|> n_0$ with $n_0$ sufficiently large 
\begin{align}\nonumber
|z\sigma_k^n|, |z\lambda_k^\pm|, |z\tau_k|, 2|z\gamma_k| \leq & 2 / 3 \qquad \forall |k|\leq |n| / 2, \: \forall |z|\leq|\tau_n|^{-1} \\\label{2.13bis}
\Big|\frac{\sigma_k^n}{\tau_k-\tau_n}\Big|, \Big|\frac{\gamma_k}{\tau_k-\tau_n}\Big|\leq& 1 /2
\qquad \forall |k|> |n| / 2, \: k\neq n \,.
\end{align}
Note that $n_0$ can be chosen uniformly on bounded subsets of $W\cap H^1_c$. 
Thus for any $|n|> n_0$ and $|z|\leq|\tau_n|^{-1},$ $g_n$ is well defined and one has 
\[
g_n(z)= \prod_{|k|\leq |n| / 2}\frac{1+ z\sigma_k^n}{\sqrt[+]{(1+z\lambda_k^+)(1+z \lambda_k^-)}}=
\prod_{|k|\leq |n| / 2}\frac{1+ \frac{\alpha_k^n}{1+z\tau_k}}{\sqrt[+]{1- \big(\frac{z\gamma_k / 2}{1+z\tau_k}\big)^2}} \, .
\]
Writing the product $\prod_{|k|\leq |n| / 2}\left(\cdot\right)$ as 
$\exp\big(\sum_{|k|\leq |n| / 2}\log(\cdot)\big)$ it follows from \cite{GK}, Lemma C.2, that 
$\sup_{|z|\leq |\tau_n|^{-1}} \left|g_n(z)-1\right|\leq C$ where $C>0$ can be chosen uniformly on bounded subsets of $W\cap H^1_c$.
\begin{lem1}\label{Lemma2.9}
For any $|n|>n_0$ with $n_0$ as in \eqref{2.13bis} and $\varphi\in W\cap H^1_c$ the coefficients of the expansion of $g_n$ at $z=0$ of order $3$, 
$g_n(z)=1+g'(0)+ g_n''(0)\frac{z^2}{2}+ r_n(z)z^3$,
are given by
$g'_n(0)=  \sum_{|k|\leq \frac{|n|}{2}} (\sigma_k^n-\tau_k)\, ,$
\[ 
g_n''(0)= g'(0)^2+ \sum_{|k|\leq |n| / 2} \big((\tau_k-\sigma_k^n)(\tau_k+\sigma_k^n)
+\big(\gamma_k / 2\big)^2\big),
\]
and $r_n(z)= \frac{1}{2}\int_0^1 g_n'''(sz)(1-s)^2ds$
where $\sup_{|z|\leq|\tau_k|^{-1}}|g_n'''(z)|\leq C$ and $C$ can be chosen uniformly on bounded subsets of $W\cap H^1_c$.
\end{lem1}
\begin{proof}[Proof of Lemma \ref{Lemma2.9}] The derivatives $g'_n(0)$ and $g''_n(0)$ are readily computed using the product rule. To estimate $g'''_n(z),$ write $g'_n(z)$ as $g_n'(z)= g_n(z)h_n(z)$ where by the product rule
\[
h_n(z)= \sum_{|k|\leq |n| / 2}
\Big(\frac{\sigma_k^n}{1+ z\sigma_k^n}- \frac{\tau_k+ \lambda_k^+\lambda_k^- z}{(1+ z\lambda_k^+)(1+ z \lambda_k^-)}\Big) \, .
\]
Then $g_n''= (g_nh_n)'= g_nh_n^2+ g_nh'_n$ whereas
$g'''_n= g_n h_n^3+ 3 g_n h_nh_n'+ g_n h_n''. $
To estimate $h_n$ and its derivatives note that
\[ 
\frac{\sigma^n_k}{1+ z\sigma^n_k}= \frac{\tau_k}{1+ z\tau_k} + \frac{\sigma^n_k-\tau_k}{(1+ z\sigma^n_k)(1+ z\tau_k)}
\]
and 
\[ 
\frac{\tau_k+ \lambda_k^+\lambda_k^- z}{(1+ z \lambda_k^+)(1+ z \lambda_k^-)}
= \frac{\tau_k}{1+ z\tau_k}- \frac{1}{1+ z\tau_k}\frac{z\left(\gamma_k / 2\right)^2}{(1+ z \lambda_k^+)(1+ z \lambda_k^-)} \, ,
\]
implying that
\[
h_n(z)=\sum_{|k|\leq |n| / 2}\Big(\frac{\alpha_k^n}{(1+ z\sigma^n_k)(1+ z\tau_k)} + 
\frac{1}{1+ z\tau_k}\frac{z\big(\gamma_k / 2\big)^2}{(1+ z \lambda_k^+)(1+ z \lambda_k^-)}\Big) \, .
\]
By our choice of $n_0$ and the asymptotics of Theorem \ref{Theorem2.4} it then follows that 
\begin{align*}
\sup_{|z|\leq|\tau_k|^{-1}}|h_n(z)|\leq & 
C_1 \sum_{|k| \leq |n| / 2} \Big(|\alpha_k^n|+|\gamma_k|^2\frac{1}{n}\Big)\leq C 
\\
\sup_{|z|\leq|\tau_k|^{-1}}|h_n'(z)|\leq &
 C_1 \sum_{|k|\leq |n| / 2} \big(\langle k\rangle|\alpha_k^n|+|\gamma_k|^2\big)\leq C
\\
\sup_{|z|\leq|\tau_k|^{-1}}|h_n''(z)|\leq & 
C_1 \sum_{|k|\leq |n| / 2} \big(\langle k\rangle^2|\alpha_k^n|+\langle k\rangle|\gamma_k|^2\big)\leq C
\end{align*}
where $C_1$ and $C$ can be chosen uniformly on bounded subsets of $W\cap H^1_c$. Together with the bound for $g_n$ established ahead of Lemma \ref{Lemma2.9}, the claimed boundedness of $\sup_{|z|\leq|\tau_k|^{-1}}|g_n'''(z)|$ follows. 
\end{proof}
Lemma \ref{Lemma2.8} and Lemma \ref{Lemma2.9} now allow us to prove the following result which completes the proof of Theorem \ref{Theorem2.7}.
\begin{lem1}\label{Lemma2.10}
For $|n|>n_0$  with $n_0$ as in \eqref{2.13bis} and $\varphi$ in $W\cap H^2_c$
\[
\sum_{|k|\leq |n| / 2}\Big( (\tau_k-\sigma_k^n)(\tau_k+\sigma_k^n)+ \Big(\frac{\gamma_k}{2}\Big)^2\Big)
=\int_0^1 \varphi_1(x)\varphi_2(x)\,dx +O\Big(\frac{1}{n}\Big)
\]
where the error term is uniform on bounded subsets of $W\cap H^2_c.$
\end{lem1}
\begin{proof}[Proof of Lemma \ref{Lemma2.10}] By the definition of $f_n$ and $g_n$ one has 
\begin{align}\label{2.14} 
f_n(0)= \prod_{k\neq n,|k|> |n| / 2 } \frac{1-\frac{\alpha^n_k}{\tau_k-\tau_n}}{\sqrt[+]{1-\left(\frac{\gamma_k / 2}{\tau_k-\tau_n}\right)^2}}
\cdot g_n\Big(-\frac{1}{\tau_n}\Big), \qquad \alpha_k^n= \sigma_k^n-\tau_k \, .
\end{align}
By Theorem \ref{Theorem2.4} and \cite{GK}, Lemma C.2, 
\begin{align*}
\Big|&\exp\Big( \sum_{k\neq n,|k|> |n| / 2}\log\Big(1+\frac{\alpha^n_k}{\tau_k-\tau_n}\Big)-1\Big)\Big|\\
\leq & \,\,\, C_1 \sum_{k\neq n,|k|> |n| / 2}\Big|\frac{\alpha^n_k}{\tau_k-\tau_n}\Big|
\leq \,\, C \, \sum_{k\neq n,|k|> |n| / 2 }|\gamma_k|^2
\end{align*}
and similarly 
\[
\Big|\exp\Big(\sum_{k\neq n,|k|> |n| / 2 }-\frac{1}{2}\log\Big(1-\Big(\frac{\frac{\gamma_k}{2}}{\tau_k-\tau_n}\Big)^2 \Big)-1\Big)\Big|
\leq \,\, C \sum_{k\neq n,|k|> |n| / 2}|\gamma_k|^2
\]
where $C_1$, $C$ can be chosen uniformly on bounded subsets of $W\cap H^2_c.$ In \cite{KST2} it is shown that $\gamma_k^2= \frac{1}{k^4}\ell^1(k)$ uniformly on bounded subsets of $W\cap H^2_c.$ Hence
\begin{align}\label{2.15} 
\prod_{k\neq n,|k|> |n| / 2 } 
\frac{1+\frac{\alpha^n_k}{\tau_k-\tau_n}}{\sqrt[+]{1-\left(\frac{\gamma_k / 2}{\tau_k-\tau_n}\right)^2}}= 
1+O\Big(\frac{1}{n^4}\Big)
\end{align}
uniformly on bounded subsets of $W\cap H^2_c.$ Substituting the asymptotic estimate \eqref{2.15} and the ones of Lemma \ref{Lemma2.8} and Lemma \ref{Lemma2.9} into the identity \eqref{2.14} one obtains
\begin{align}\label{2.16}
1+ O\Big( \frac{1}{n^3}\Big) = \Big(1+ O\Big( \frac{1}{n^4}\Big)\Big)\cdot 
\left(1+\sum_{|k|\leq |n| / 2}(\sigma_k^n-\tau_k)\Big(-\frac{1}{\tau_n}\Big)+ A 
+ O\Big( \frac{1}{n^3}\Big)\right)
\end{align}
uniformly on bounded subsets of $W\cap H^2_c$  where
\[
A := \frac{1}{2\tau_n^2} \Big(\sum_{|k|\leq |n| / 2}(\sigma_k^n-\tau_k)\Big)^2   + 
\frac{1}{2\tau_n^2}\Big(\sum_{|k|\leq |n| / 2}(\tau_k-\sigma_k^n)(\tau_k+\sigma_k^n) 
+\big(\gamma_k / 2\big)^2 \Big).
\]
By  \cite{GK1}, Lemma 5.8, one knows that
$\sum_{k\neq n} (\sigma_k^n-\tau_k) = \tau_n -n\pi$. As 
\[
\sum_{k|> |n| / 2}|\sigma_k^n-\tau_k|\leq C \sum_{k|> |n| / 2} |\gamma_k|^2
= O\Big(\frac{1}{n^2}\Big)
\] 
uniformly on bounded subsets of $W\cap H^1_c$ one has 
\[
\sum_{k|\leq |n| / 2}(\sigma_k^n-\tau_k)\Big(-\frac{1}{\tau_n}\Big) = 
(\tau_n-n\pi)\Big(-\frac{1}{\tau_n}\Big) +O\Big(\frac{1}{n^3}\Big) .
\]
As by Theorem \ref{Theorem2.4}, 
$\tau_n-n\pi= \frac{1}{2\pi n}\int_0^1 \varphi_1(x)\varphi_2(x)\,dx + O\left(\frac{1}{n^2}\right)$ it 
then follows 
\[
\sum_{k|\leq |n| / 2}(\sigma_k^n-\tau_k)\Big(-\frac{1}{\tau_n}\Big) = 
-\frac{1}{2\pi n\tau_n}\int_0^1 \varphi_1(x)\varphi_2(x) dx + O\Big(\frac{1}{n^3}\Big)
\]
whereas 
$\big(\sum_{k|\leq |n| / 2}(\sigma_k^n-\tau_k)(-\frac{1}{\tau_n})\big)^2 = O\big(\frac{1}{n^4}\big).$
When substituted into \eqref{2.16}, these estimates lead to 
\[
-\int_0^1 \varphi_1(x)\varphi_2(x) dx+\sum_{|k|\leq |n| /2}(\tau_k-\sigma_k^n)(\tau_k+\sigma_k^n)+
\Big(\frac{\gamma_k}{2}\Big)^2= O\Big(\frac{1}{n}\Big)
\] 
uniformly on bounded subsets of $W\cap H^2_c$ as claimed.
\end{proof}

\section{Proof of the main results}
In this section we prove Theorem \ref{Theorem1.1}, Corollary \ref{Corollary1.2}, and Theorem \ref{Theorem1.3} as well as some additional results. Without further reference we use the notation and terminology introduced in the previous sections. We begin with the proof of Theorem \ref{Theorem1.1}. 
First  we review the results on the semi-linearity of the Birkhoff map $\Phi$ relative to the Fourier transform $\mathcal{F}$, established in \cite{KSTnls}. We say that a (possibly non-linear) map between (subsets of) Banach spaces is bounded if the image of any bounded subsets (in the domain of definition) is bounded.
\begin{thm1}\label{Theorem3.1} \cite{KSTnls} 
For any $N\in \mathbb{Z}_{\geq 1},$ the restriction of $\Phi-\mathcal{F}$ to $H^N_r$ takes values in $\mathfrak{h}^{N+1}_r.$ The map $\Phi-\mathcal{F}: H^{N}_r\to \mathfrak{h}^{N+1}_r$ is real analytic and bounded.
\end{thm1}
\noindent Theorem \ref{Theorem3.1} leads to the following two corollaries:
\begin{cor1}\label{Corollary3.2} \cite{KSTnls} (i)
For any $N\in \mathbb{Z}_{\geq 1},$ the restriction of $\Phi^{-1}- \mathcal{F}^{-1}$ to $\mathfrak{h}^N_r$ takes values in $H^{N+1}_r$. The map  $\Phi^{-1}- \mathcal{F}^{-1}: \mathfrak{h}^N_r\to H^{N+1}_r$ is real analytic and bounded.
(ii) For any
 $s \in \mathbb{R}_{\geq 1},$ the restriction of 
 $\Phi$ to $H^s_r$ takes values in  $\mathfrak{h}^s_r$. 
The map $\Phi:H^s_r\to  \mathfrak{h}^s_r$ is a real analytic and bounded diffeomorphism as is its inverse.
\end{cor1}
\begin{cor1}\label{Corollary3.3} \cite{KSTnls}
For any $s\in \mathbb{R}_{\geq 1}, \, \Phi:H^s_r\to \mathfrak{h}^s_r$ and $\Phi^{-1}: \mathfrak{h}^s_r \to H^s_r$ are weakly continuous. 
\end{cor1}


\noindent Denote by $S^t, \, t\in \mathbb{R},$ the dNLS flow on $H^N_r,$ i.e., 
$S^t:H^N_r\to H^N_r, \,\varphi \mapsto S^t(\varphi)$ where 
$t\mapsto S^t(\varphi)$
denotes the solution of \eqref{1}-\eqref{2} for the initial data $\varphi= (u,\bar u)\in H^N_r.$ 
Furthermore define 
\[
R^t:H^N_r\to H^N_r,\; \varphi \mapsto S^t(\varphi)- \big(v(t), \bar v(t)\big)
\] 
where $v(x,t)=\sum_{n\in \mathbb{Z}} \hat \varphi_1(n)e^{- i\omega_n^{\rm NLS}t} e^{2\pi inx} $ and $\omega_n^{\rm NLS}= \omega_n^{\rm NLS}(\varphi).$
\begin{thm1}\label{Theorem3.4}
For any $N\in \mathbb{Z}_{\geq 1}$ and any $t\in \mathbb{R},$ the flow $R^t$ maps $H_r^N$ into $H^{N+1}_r.$ It has the following properties:
\begin{itemize}
\item[(i)] for any $\varphi\in H^N_r,$ the orbit $\left\{R^t\left(\varphi\right) |\, t \in \mathbb{R}\right\}$ is relatively compact in $H_r^{N+1};$
\item[(ii)]for any $M>0$, the union of orbits $\left\{R^t\left(\varphi\right) |\, t \in \mathbb{R},\, \varphi\in H^N_r,\, \|\varphi\|_{H^N}\leq M \right\}$ is bounded in $H_r^{N+1};$
\item[(iii)] $R^t:H^N_r\to H^{N+1}_r$ is real analytic and bounded;
\item[(iv)] for $\varphi \in H^N_r, \; \partial_t R^t(\varphi)$ takes values in $H_r^{N-1}$ and the map 
$\partial_t R^t:H^N_r\to H^{N-1}_r$ is real analytic and bounded.
\end{itemize}
\end{thm1}
\begin{proof}[Proof of Theorem \ref{Theorem3.4}] 
Fix $N\in \mathbb{Z}_{\geq 1}$ and $t\in \mathbb{R}.$ In Section 2 we have introduced the real analytic map
$\Omega^t: H^1_r\to \ell_r^\infty.$ Note that $\Omega^t$ is related to the flow map in Birkhoff coordinates as follows
\begin{align}
\label{3.1} S^t= \Phi^{-1}\circ \Omega^t\cdot \Phi
\end{align}
where $\Omega^t\cdot \Phi$ is the map
$H^N_r \to \mathfrak{h}_r^N, \varphi \mapsto \Omega^t(\varphi)\cdot \Phi(\varphi) $
given by
\[
\Omega^t(\varphi)\cdot \Phi(\varphi) =
\Big(\big(\Omega^t_1(n,\varphi)z_1(n,\varphi)\big)_{n\in \mathbb{Z}}, 
\big(\Omega^t_2(n,\varphi)z_2(n,\varphi)\big)_{n\in \mathbb{Z}}\Big)\,.
\]
Clearly, the map 
\[
\ell^\infty_c\times \mathfrak{h}_c^s\to \mathfrak{h}_c^s, \,\, 
\big((w_1,w_2),(z_1,z_2)\big) \mapsto 
\big((w_1(n) \cdot z_1(n))_{n\in \mathbb{Z}},(w_2(n) \cdot z_2(n))_{n\in \mathbb{Z}}\big) 
\]
is bilinear and real analytic.
Write $\Phi$ and $\Phi^{-1}$ in the form $\Phi= \mathcal{F}+ A$ respectively $\Phi^{-1}= \mathcal{F}^{-1}+ B$ and substitute 
these expressions into formula \eqref{3.1} for $S^t$ to get 
$S^t=\left(\mathcal{F}^{-1}+ B\right)\circ \Omega^t \cdot \left(\mathcal{F}+ A\right)
= \mathcal{F}^{-1}\circ \Omega^t \cdot \mathcal{F} + R^t $
where \begin{align}
\label{3.2} R^t= \mathcal{F}^{-1}\circ \Omega^t \cdot  A +B\circ \Omega^t \cdot \Phi.
\end{align}
By Theorem \ref{Theorem3.1} and Corollary \ref{Corollary3.2} (i),  $A:H^N_r \to \mathfrak{h}^N_r$ 
respectively $B:\mathfrak{h}^N_r \to H^{N+1}_r$ are real analytic and bounded. It then follows that
\[ 
\Omega^t\cdot \Phi: H^N_r \to \mathfrak{h}_r^N, \, 
\varphi \mapsto \Omega^t(\varphi)\cdot \Phi(\varphi)\, ,
\quad
 \Omega^t\cdot A: H^N_r \to \mathfrak{h}_r^{N+1}, \, 
\varphi \mapsto \Omega^t(\varphi)\cdot A(\varphi)  
\]
are both real analytic and bounded. Hence $R^t$ is the sum of real analytic bounded maps which proves item (iii). As $\|\Omega^t(\varphi)\|_{\ell^\infty}=1$ for any $\varphi\in H^N_r,t\in \mathbb{R}$ and $\Phi:H^N_r \to \mathfrak{h}_r^N$ as well as $B:\mathfrak{h}_r^N \to H^{N+1}_r$ are bounded, item (ii)
 is proved as well. 
 In view of the characterisation of relatively compact subsets of $\mathfrak{h}^s_r$, item (i) also follows.
Towards (iv), note that by \eqref{3.2}, $\partial_tR^t$ can be computed as
\[
\partial_tR^t(\varphi)= \mathcal{F}^{-1}\circ \partial_t \Omega^t\cdot A\big|_\varphi + 
dB \big |_{(\Omega^t\cdot \Phi)(\varphi) }\circ \left(\partial_t \Omega^t\cdot \Phi\right)(\varphi).
\]
Clearly, 
$\partial_t \Omega^t(\varphi) = 
\Big(\big(-i\omega_n^{\rm NLS}(\varphi)\Omega_1(n,\varphi)\big)_{n\in \mathbb{Z}},
\big(i\omega_n^{\rm NLS}(\varphi)\Omega_2(n,\varphi)\big)_{n\in \mathbb{Z}}\Big).$
Arguing as in the proof of items (ii) and (iii) and taking into account the asymptotics $\omega_n^{\rm NLS}= 4\pi^2n^2+ O(1)$ as $n\to \infty$ of Corollary \ref{Corollary2.5}, it follows that the map $\partial_t R^t:H^N_r\to H^{N-1}_r$ possesses the stated properties.  
 \end{proof}
 \begin{proof}
 [Proof of Theorem \ref{Theorem1.1}] The claimed statement is contained in Theorem \ref{Theorem3.4}.
 \end{proof}
 \begin{proof}[Proof of Corollary \ref{Corollary1.2}]
 Let $N=\lfloor s\rfloor$. Using the estimates
 \[
 \|u(\cdot,t)\|_{H^s}\leq \|v(\cdot,t)\|_{H^s}+ \|u(\cdot,t)-v(\cdot,t)\|_{H^s},\,  
 \|v(\cdot,t)\|_{H^s} =  \|v(\cdot,0)\|_{H^s}=\!\|u_0\|_{H^s}
 \]
 and $\|u(\cdot,t)-v(\cdot,t)\|_{H^s} \leq \|u(\cdot,t)-v(\cdot,t)\|_{H^{N+1}} $
 it then follows from Theorem \ref{Theorem1.1} that $\sup_{t\in \mathbb{R}} \|u(\cdot,t)\|_{H^s}\leq C$ where $C>0$ can be chosen 
 uniformly for bounded sets of initial data $u_0$ in $H^s$.
 \end{proof}
 \begin{proof}
 [Proof of Theorem \ref{Theorem1.3}] Let $u_0\in H^N_r$ with $N\geq 2$. Then 
 \[\|u(\cdot,t)-w(\cdot,t)\|_{H^{N+1}} \leq \|u(\cdot,t)-v(\cdot,t)\|_{H^{N+1}} + \|v(\cdot,t)-w(\cdot,t)\|_{H^{N+1}} \]
 By Theorem \ref{Theorem1.1}
 \begin{align}
 \label{3.3} \sup_{t\in \mathbb{R}}\|u(\cdot,t)-v(\cdot,t)\|_{H^{N+1}}\leq C_1 
 \end{align}
 where $C_1>0$ can be chosen uniformly on bounded sets of initial data in $H^N_r.$
 Furthermore, with $\tilde\omega_n := 4\pi^2 n^2 + 4 \int_0^1 |u_0(x)|^2 dx,$ one has 
 \[\|v(\cdot,t)-w(\cdot,t)\|_{H^{N+1}}^2= \sum_{n\in \mathbb{Z}} \langle n\rangle^{2N+2}|\hat u_0(n)|^2 
\big|e^{i(\omega_n^{\rm NLS}-\tilde\omega_n)t}-1\big|^2.\]
 As 
\[
\big|e^{i(\omega_n^{\rm NLS}-\tilde\omega_n)t}-1\big|= 
\Big|i (\omega_n^{\rm NLS}-\tilde\omega_n)\int_0^te^{i(\omega_n^{\rm NLS}-\tilde\omega_n)y}\,dy\Big|
  \leq \big|\omega_n^{\rm NLS}-\tilde\omega_n\big||t|
\]
it follows from Theorem \ref{Theorem1.4} that 
\begin{align}
\label{3.4} \|v(\cdot,t)-w(\cdot,t)\|_{H^{N+1}} \leq C_2 |t| \|u_0\|_{H^N}
\end{align} 
where $C_2>0$ can be chosen uniformly on bounded sets of initial data in $H^N_r.$ Combining \eqref{3.3} and \eqref{3.4} leads to the claimed estimate. 
 \end{proof}
\noindent We finish this section with an additional application of the results obtained so far. It says that the solution map $S^t$ of 
\eqref{1}-\eqref{2} is weakly continuous:
\begin{cor1}\label{Corollary3.5}
For any $s\in \mathbb{R}_{\geq 2},$ the flow map $S^t:H^s_r \to H^s_r$ of \eqref{1}-\eqref{2} is weakly continuous.
\end{cor1} 
\begin{rem1}\label{Remark3.6}
By a refined analysis of the asymptotics of the dNLS frequencies, Corollary \ref{Corollary3.5} can be shown to hold for $s\geq 0.$
Recently, for $s=0,$ weak continuity of $S^t$ has been established in \cite{OS} by very different methods than the ones used in this paper. 
\end{rem1}
\begin{proof}
[Proof of Corollary \ref{Corollary3.5}] Given $s\in \mathbb{R}_{\geq 2}$ let $N=\lfloor s\rfloor \in \mathbb{Z}_{\geq 2}$. 
Assume that the sequence $\left(\varphi^{(j)}\right)_{j\geq 1}$ is in $H^s_r$, converging weakly to an element $\varphi \in H^s_r.$  Let $z^{(j)}:=\Phi(\varphi^{(j)})$ and $z:= \Phi(\varphi).$
By Corollary \ref{Corollary3.3}, $z^{(j)} \rightharpoondown z$ converges weakly in $\mathfrak{h}^s_r.$ In particular 
\[z_1^{(j)}(n)\underset{j\to \infty}{\longrightarrow} z_1(n)\quad \text{and} \quad z_2^{(j)}(n) \underset{j\to \infty}{\longrightarrow} z_2(n)\quad \forall n\in \mathbb{Z}.\]
To see that 
\begin{align}
\label{3.5} S^t\big(\varphi^{(j)}\big)= \Phi^{-1}\big(\Omega^t\big(\varphi^{(j)}\big)\cdot z^{(j)}\big)
\end{align} 
converges
weakly in $H^s_r$ to $S^t(\varphi),$ note that by Rellich's theorem and the assumption $s\geq 2$ one has 
$\varphi^{(j)} \underset{j\to \infty}{\longrightarrow} \varphi$ strongly in $H^1_r$ implying that for any $n\in \mathbb{Z},\; \omega_n^{\rm NLS}\left(\varphi^{(j)}\right) \underset{j\to \infty}{\longrightarrow}\omega_n^{\rm NLS}( \varphi)$ by Corollary \ref{Corollary2.5}. It then follows that for any
$n\in \mathbb{Z},$
\[\lim_{j\to \infty} e^{-i \omega_n^{\rm NLS}(\varphi^{(j)})t}z_1^{(j)}(n)= e^{-i \omega_n^{\rm NLS}(\varphi)t}z_1^{(j)}(n)\]
and 
\[\lim_{j\to \infty} e^{i \omega_n^{\rm NLS}(\varphi^{(j)})t}z_2^{(j)}(n)= e^{i \omega_n^{\rm NLS}(\varphi)t}z_2^{(j)}(n).\]
As $\|\Omega^t(\varphi^{(j)})\cdot z^{(j)} \|_s=\|z^{(j)} \|_s\leq C $ for some constant $C>0$ independent of $j\geq 1$ it then follows that 
\[\Omega^t(\varphi^{(j)})\cdot z^{(j)}\underset{j\to \infty}{\rightharpoondown} \Omega^t(\varphi)\cdot z \quad \text{weakly in} \; \mathfrak{h}^s_r.\]
Using that by Corollary \ref{Corollary3.3}, $\Phi^{-1}:\mathfrak{h}^s_r\to H^s_r$ is weakly continuous, $S^t(\varphi^{(j)})\rightharpoondown S^t(\varphi)$ by \eqref{3.5}.
\end{proof}


\end{document}